\documentclass[leqno,12pt]{article}

\usepackage{latexsym}
 \usepackage[dvips]{graphicx}
\usepackage{amsmath}
\usepackage{amssymb}
\usepackage{cite}

\newcommand{\nc}{\newcommand}
\nc{\nt}{\newtheorem}
\nt{thm}{Theorem}[section]
\nt{cor}[thm]{Corollary}
\nt{prop}[thm]{Proposition}
\nt{lem}[thm]{Lemma}
\nt{defn}[thm]{Definition}
\nt{rem}[thm]{Remark}
\nt{exa}[thm]{Example}
\nt{ass}[thm]{Assumption}
\nc{\ip}[2]{\mbox{$\langle #1,#2 \rangle$}}
\nc{\pf}{\noindent{\bf Proof\ \ }}
\nc{\finpf}{\hfill{$\Box$}\linespace}
\nc{\linespace}{\vspace{\baselineskip} \noindent}
\nc{\R}{{\bf R}}
\nc{\C}{{\bf C}}
\nc{\E}{{\bf E}}
\nc{\Y}{{\bf Y}}
\nc{\RE}{\mbox{\rm Re}\,}
\nc{\Rn}{{\bf R}^n}
\nc{\Mn}{{\bf M}^n}
\nc{\bx}{\bar{x}}
\nc{\by}{\bar{y}}
\nc{\e}{\epsilon}
\nc{\inT}{\mbox{\rm int}\,}
\nc{\cl}{\mbox{\rm cl}\,}
\nc{\gph}{\mbox{\rm gph}\,}
\nc{\conv}{\mbox{\rm conv}\,}
\nc{\rt}{\rightarrow}
\nc{\tra}{\mbox{\rm tr}\,}

\nc{\xbar}{{\overline{x}}}
\nc{\vbar}{{\overline{v}}}
\nc{\yhat}{{\widehat{y}}}
\nc{\xhat}{{\widehat{x}}}

\nc{\und}{\quad\mbox{ and }\quad}
\nc{\emp}{\ensuremath{\varnothing}}

\def\tto{\;{\lower 1pt \hbox{$\rightarrow$}}\kern -12pt
           \hbox{\raise 2.8pt \hbox{$\rightarrow$}}\;}
\newenvironment{myequation}{\setcounter{equation}{\value{thm}}
   \begin{equation}}{\addtocounter{thm}{1}\end{equation}}
\newenvironment{myeqnarray}{\setcounter{equation}{\value{thm}}
    \begin{eqnarray}}{\setcounter{thm}{\value{equation}}\end{eqnarray}}

\nc{\bmye}{\begin{myequation}}
\nc{\emye}{\end{myequation}}

\newcommand{\trans}[1]{{#1}^{\top}}

\newcommand{\norm}[1]{\|#1\|}
\newcommand{\norminf}[1]{\norm{#1}_{\infty}}

\begin{document}
\title{
Local convergence for alternating and averaged nonconvex projections}
\author{
 A.S. Lewis\thanks{ORIE, Cornell University, Ithaca, NY 14853, U.S.A.
\texttt{aslewis\char64 orie.cornell.edu people.orie.cornell.edu/\~~\hspace{-4pt}aslewis}.
Research supported in part by National Science Foundation Grant DMS-0504032.}
\and
D.R. Luke\thanks{
Department of Mathematical Sciences, University of Delaware. \hfill \mbox{} \texttt{rluke@math.udel.edu}}
\and
J. Malick\thanks{CNRS, Lab. Jean Kunztmann, University of Grenoble. \texttt{jerome.malick@inria.fr}}
}
\maketitle


\medskip


\begin{abstract}
The idea of a finite collection of closed sets having ``strongly regular intersection'' at a given point is crucial in variational analysis.  We show that this central theoretical tool also has striking algorithmic consequences.  Specifically, we consider the case of two sets, one of which we assume to be suitably ``regular'' (special cases being convex sets, smooth manifolds, or feasible regions satisfying the Mangasarian-Fromovitz constraint qualification).  We then prove that von Neumann's method of ``alternating projections'' converges locally to a point in the intersection, at a linear rate associated with a modulus of regularity.   As a consequence, in the case of several arbitrary closed sets having strongly regular intersection at some point, the method of ``averaged projections'' converges locally at a linear rate to a point in the intersection.  Inexact versions of both algorithms also converge linearly.
\end{abstract}

\section{Introduction}

An important theme in computational mathematics is the relationship between the ``conditioning'' of a problem instance and the speed of convergence of iterative solution algorithms on that instance.  A classical example is the method of conjugate gradients for solving a positive definite system of linear equations:  we can bound the linear convergence rate in terms of the relative condition number of the associated matrix.  More generally, Renegar \cite{Ren95,Ren95a,Ren96} showed that the rate of convergence of interior-point methods for conic convex programming can be bounded in terms of the ``distance to ill-posedness'' of the program.

In studying the convergence of iterative algorithms for nonconvex minimization problems or nonmonotone variational inequalities, we must content ourselves with a local theory.  A suitable analogue of the distance to ill-posedness is then the notion of ``metric regularity'', fundamental in variational analysis.  Loosely speaking, a generalized equation, such as a system of inequalities, for example, is metrically regular when, locally, we can bound the distance from a trial solution to an exact solution by a constant multiple of the error in the equation generated by the trial solution.  The constant needed is called the ``regularity modulus'', and its reciprocal has a natural interpretation as a distance to ill-posedness for the equation \cite{Don03}.

This philosophy suggests understanding the speed of convergence of algorithms for solving generalized equations in terms of the regularity modulus at a solution.  Recent literature focuses in particular on the proximal point algorithm (see for example \cite{Pen02,Ius03,Ara05}).  A unified approach to the relationship between metric regularity and the linear convergence of a family of conceptual algorithms appears in \cite{Kla07}.

We here study a very basic algorithm for a very basic problem.  We consider the problem of finding a point in the intersection of several closed sets, using the {\em method of averaged projections}:  at each step, we project the current iterate onto each set, and average the results to obtain the next iterate.  Global convergence of this method in the case of two closed convex sets was proved in 1969 in \cite{Aus69}.  In this work we show, in complete generality, that this method converges locally to a point in the intersection of the sets, at a linear rate governed by an associated regularity modulus.  Our linear convergence proof is elementary:  although we use the idea of the normal cone, we apply only the definition, and we discuss metric regularity only to illuminate the rate of convergence.

Our approach to the convergence of the method of averaged projections is standard \cite{Pie84,Bau93}: we identify the method with von Neumann's {\em alternating} projection algorithm \cite{von50} on two closed sets (one of which is a linear subspace) in a suitable product space.  A nice development of the classical method of alternating projections may be found in \cite{Deu01}.  The linear convergence of the method for two closed convex sets with regular intersection was proved in \cite{Bau93}, strengthening a classical result of \cite{Gub67}.  Remarkably, we show that, assuming strong regularity, {\em local} linear convergence requires good geometric properties (such as convexity, smoothness, or more generally, ``amenability '' or ``prox-regularity'') of only one of the two sets.

One consequence of our convergence proof is an algorithmic demonstration of the ``exact extremal principle'' described in \cite[Theorem 2.8]{Mor06}.  This result, a unifying theme in \cite{Mor06}, asserts that if several sets have strongly regular intersection at a point, then that point is not ``locally extremal'' \cite{Mor06}:  in other words, translating the sets by small vectors cannot render the intersection empty locally.  To prove this result, we simply apply the method of averaged projections, starting from the point of regular intersection.  In a further section, we show that inexact versions of the method of averaged projections, closer to practical implementations, also converge linearly.

The method of averaged projections is a conceptual algorithm that might appear hard to implement on concrete nonconvex problems.  However, the projection problem for some nonconvex sets is relatively easy.  A good example is the set of matrices of some fixed rank:  given a singular value decomposition of a matrix, projecting it onto this set is immediate.  Furthermore, nonconvex iterated projection algorithms and analogous heuristics are quite popular in practice, in areas such as inverse eigenvalue problems \cite{chu-1995,chen-chu-1996}, pole placement \cite{orsi-2006,yang-orsi-2006}, information theory \cite{tropp-dhillon-heath-strohmer-2005}, low-order control design \cite{grigoriadis-skelton-1996,grigoriadis-beran-2000,orsi-helmke-moore-2005} and image processing \cite{weber-allebach-1986,bauschke-combettes-luke-2002}).  Previous convergence results on nonconvex alternating projection algorithms have been uncommon, and have either focussed on a very special case (see for example \cite{chen-chu-1996,Lew06Alternating}), or have been much weaker than for the convex case \cite{combettes-trussell-1990,tropp-dhillon-heath-strohmer-2005}.  
For more discussion, see \cite{Lew06Alternating}.

Our results primarily concern R-linear convergence:  in other words, we show that our sequences of iterates converge, with error bounded by a geometric sequence.  In a final section, we employ a completely different approach to show that the method of averaged projections, for prox-regular sets with regular intersection, has a Q-linear convergence property:  each iteration guarantees a fixed rate of improvement.  In a final section, we illustrate these theoretical results with an elementary numerical example coming from signal processing.


\section{Notation and definitions}

We begin by fixing some notation and definitions.  Our underlying setting throughout this work is a Euclidean space $\E$ with corresponding closed unit ball $B$.  For any point $x \in \E$ and radius $\rho > 0$ , we write $B_{\rho}(x)$ for the set $x+\rho B$.

Consider first two sets $F,G \subset \R$.  A point $\bx \in F \cap G$ is {\em locally extremal} \cite{Mor06} for this pair of sets if restricting to a neighborhood of $\bx$ and then translating the sets by small distances can render their intersection empty:  in other words, there exists a $\rho > 0$ and a sequence of vectors $z_r \rt 0$ in $\E$ such that
\[
(F + z_r) \cap G \cap B_{\rho}(\bx) ~=~ \emp ~~\mbox{for all}~r=1,2,\ldots .
\]
Clearly $\bx$ is not locally extremal if and only if
\[
0 \in \inT \Big( ((F  - \bx) \cap \rho B) - ((G - \bx) \cap \rho B) \Big)
~~\mbox{for all}~ \rho > 0.
\]

For recognition purposes, it is easier to study a weaker property than local extremality.
Following the terminology of  \cite{Kru06}, we say
{\em the two sets $F,G \subset \E$ have  strongly regular intersection} at the point
$\bx \in F \cap G$ if there exists a constant $\alpha > 0$ such that
\[
\alpha \rho B ~\subset~ ((F-x) \cap \rho B) - ((G-z) \cap \rho B)
\]
for all points $x \in F$ near $\bx$ and $z \in G$ near $\bx$.  By considering the case  $x=z=\bx$, we see that strong regularity implies that $\bx$ is not locally extremal.  This ``primal'' definition of strong regularity is often not the most convenient way to handle strong regularity, either conceptually or theoretically.  By contrast, a ``dual'' approach, using normal cones, is very helpful.

Given a set $F \subset \E$, we define the {\em distance function} and (multivalued) {\em projection} for $F$ by
\begin{eqnarray*}
d_F(x) & = & d(x,F) ~=~ \inf \{ \|z-x\| : z \in F \} \\
P_F(x) & = & \mbox{argmin} \{ \|z-x\| : z \in F \}.
\end{eqnarray*}
The central tool in variational analysis is the {\em normal cone} to a closed set
$F \subset \E$ at a point $\bx \in F$, which can be defined (see \cite{Cla98,Mor06,Roc98}) as
\[
N_F(\bx) ~=~ \Big\{ \lim_i t_i(x_i - z_i) : t_i \ge 0,~ x_i \rt \bx,~ z_i \in P_F(x_i) \Big\}.
\]
Notice two properties in particular.  First,
\bmye \label{projections-normals}
z \in P_F(x) ~~\Rightarrow~~ x-z \in N_F(z).
\emye
Secondly, the normal cone is a ``closed'' multifunction:  for any sequence of points $x_r \rt \bx$ in $F$, any limit of a sequence of normals $y_r \in N_F(x_r)$ must lie in $N_F(\bx)$.  Indeed, the definition of the normal cone is in some sense driven by these two properties:  it is the smallest cone satisfying the two properties.  Notice also that we have the equivalence: $N_F(x)=\{0\}\iff x\in \inT F$.

Normal cones provide an elegant alternative approach to defining strong regularity.
In general, {\em a family of closed sets $F_1,F_2,\ldots F_m \subset \E$
has strongly regular intersection} at a point $\bx \in \cap_i F_i$,
if the only solution to the system
\begin{eqnarray*}
y_i & \in & N_{F_i}(\bx) ~~ (i=1,2,\ldots,m) \\
\displaystyle{\sum_{i=1}^m} y_i & = & 0,
\end{eqnarray*}
is $y_i = 0$ for $i=1,2,\ldots,m$.  In the case $m=2$, this condition can be written
\[
N_{F_1}(\bar x) \cap - N_{F_2}(\bar x) =\{0\},
\]
and it is equivalent to our previous definition (see \cite[Cor 2]{Kru06}, for example).
We also note that this condition appears throughout variational-analytic theory.
For example, it guarantees the important inclusion
(see \cite[Theorem 6.42]{Roc98})
\[
N_{F_1 \cap \ldots \cap F_m}(\bx) ~\subset~ N_{F_1}(\bx) + \cdots +  N_{F_m}(\bx).
\]

We will find it helpful to quantify the notion of strong regularity (cf.\ \cite{Kru06}).  A straightforward compactness argument shows the following result.

\begin{prop} [quantifying strong regularity]
  A collection of closed sets $F_1,F_2,\ldots,F_m \subset \E$
  have strongly regular intersection at a point $\bx \in \cap F_i$
  if and only if there exists a constant $k > 0$ such that the following condition holds:
  \bmye \label{conditioning}
  y_i \in N_{F_i}(\bx) ~~(i=1,2,\ldots,m)~~ \Rightarrow~~
  \sqrt{\sum_i \|y_i\|^2} \le k \Big\| \sum_i y_i \Big\|.
  \emye
\end{prop}

\noindent
We define the {\em condition modulus} $\mbox{cond}(F_1,F_2,\ldots,F_m | \bx)$ to be the infimum of all constants $k > 0$ such that property (\ref{conditioning}) holds.  Using the triangle and Cauchy-Schwarz inequalities,
we notice that  vectors $y_1,y_2,\ldots,y_m \in \E$ always satisfy the inequality
\bmye \label{cauchy}
\sum_i\| y_i \|^2 \ge \frac{1}{m} \Big\| \sum_i y_i \Big\|^2,
\emye
which yields
\bmye \label{e:boundoncond}
\mbox{cond}(F_1,F_2,\ldots,F_m | \bx) \ge \frac{1}{\sqrt{m}}\,,
\emye
except in the special case when $N_{F_i}(\bx) = \{0\}$
(or equivalently $\bx \in \inT F_i$) for all $i=1,2,\ldots,m$;
in this case the condition modulus is zero.

 One goal of this paper is to show that,
far from being of purely analytic significance, strong regularity has central algorithmic consequences, specifically for the method of averaged projections for finding a point in the intersection $\cap_i F_i$.  Given any initial point $x_0 \in \E$, the algorithm proceeds iteratively as follows:
\begin{eqnarray*}
  z^i_n & \in & P_{F_i}(x_n) ~~~ (i=1,2,\ldots,m) \\
  x_{n+1}   & = & \frac{1}{m}(z^1_n + z^2_n + \cdots + z^m_n).
\end{eqnarray*}
Our main result shows, assuming only strong regularity, that providing the initial point $x_0$ is near $\bx$, any sequence $x_1,x_2,x_3,\ldots$ generated by the method of averaged projections converges linearly to a point in the intersection $\cap_i F_i$, at a rate governed by the condition modulus.

\section{Strong and metric regularity}

The notion of strong regularity is well-known to be closely related to another central idea in variational analysis:  ``metric regularity''.  A concise summary of the relationships between a variety of regular intersection properties and metric regularity appears in \cite{Kru06}.  We summarize the relevant ideas here.

Consider a set-valued mapping $\Phi \colon \E \tto \Y$, where $\Y$ is a second Euclidean space.  The inverse mapping $\Phi^{-1} \colon \Y \tto \E$ is defined by
\[
x \in \Phi^{-1}(y) ~\Leftrightarrow~ y \in \Phi(x),~~ \mbox{for}~ x \in \E,~y \in \Y.
\]
For vectors $\bx \in \E$ and $\bar y \in \Phi(\bx)$, we say $\Phi$ is {\em metrically regular} at $\bx$ for $\bar y$ if there exists a constant $\kappa > 0$ such that all vectors $x \in \E$ close to $\bx$ and vectors $y \in \Y$ close to $\bar y$ satisfy
\[
d(x,\Phi^{-1}(y)) \le \kappa d(y,\Phi(x)).
\]
Intuitively, this inequality gives a local linear bound for the distance to a solution of the generalized equation $y \in \Phi(x)$ (where the vector $y$ is given and we seek the unknown vector $x$), in terms of the the distance from $y$ to the set $\Phi(x)$.  The infimum of all such constants $\kappa$ is called the {\em modulus of metric regularity} of $\Phi$ at $\bx$ for $\bar y$, denoted $\mbox{reg}\,\Phi(\bx|\bar y)$.  This modulus is a measure of the sensitivity or ``conditioning'' of the generalized equation $y \in \Phi(x)$.  To take one simple example, if $\Phi$ is a single-valued linear map, the modulus of regularity is the reciprocal of its smallest singular value.  In general, variational analysis provides a powerful calculus for computing the regularity modulus.  In particular, we have the following formula \cite[Thm 9.43]{Roc98}:
\bmye \label{coderivative}
\frac{1}{\mbox{reg}\,\Phi(\bx|\bar y)} ~=~
\min \Big\{ d(0,D^*\Phi(\bx|\bar y)(w)) : w \in \Y,~ \|w\|=1 \Big\},
\emye
where $D^*$ denotes the ``coderivative''.

 We now study these ideas for a particular mapping,
highlighting the connections between metric and strong regularity.
As in the previous section, consider closed sets $F_1,F_2,\ldots,F_m \subset \E$ and a point $\bx \in \cap_i F_i$.  We endow the space $\E^m$ with the inner product
\[
\Big\langle (x_1,x_2,\ldots,x_m) , (y_1,y_2,\ldots,y_m) \Big\rangle
~=~ \sum_i \ip{x_i}{y_i},
\]
and define set-valued mapping $\Phi \colon \E \tto \E^m$ by
\[
\Phi(x) ~=~ (F_1 - x) \times (F_2 - x) \times \cdots \times (F_m - x).
\]
Then the inverse mapping is given by
\[
\Phi^{-1}(y) = \bigcap_i (F_i - y_i),~~ \mbox{for}~ y \in \E^m
\]
and finding a point in the intersection $\cap_i F_i$ is equivalent to finding a solution of the generalized equation $0 \in \Phi(x)$.  By definition, the mapping $\Phi$ is metrically regular at $\bx$ for $0$ if and only if there is a constant $\kappa > 0$ such that the following {\em strong metric inequality} holds:
\bmye \label{strong}
~~~~ d\Big(x,\bigcap_i(F_i - z_i)\Big) ~\le~ \kappa \sqrt{ \sum_i d^2(x,F_i - z_i) }
~~\mbox{for all}~ (x,z)~ \mbox{near}~ (\bx,0).
\emye
Furthermore, the regularity modulus $\mbox{reg}\,\Phi(\bx|0)$ is just the infimum of those constants $\kappa > 0$ such that inequality (\ref{strong}) holds.

To compute the coderivative $D^*\Phi(\bx|0)$, we decompose the mapping $\Phi$ as $\Psi - A$, where, for points $x \in \E$,
\begin{eqnarray*}
\Psi(x) & = & F_1 \times F_2 \times \cdots \times F_m \\
Ax & = & (x,x,\ldots,x).
\end{eqnarray*}
The calculus rule \cite[10.43]{Roc98} yields $D^*\Phi(\bx|0) = D^*\Psi(\bx|A\bx) - A^*$.
Then, by definition,
\[
v \in D^*\Psi(\bx|A\bx)(w) ~\Leftrightarrow~
(v,-w) \in N_{\mbox{\scriptsize gph}\, \Psi}(\bx,A\bx),
\]
and since $\mbox{gph}\, \Psi = \E \times F_1 \times F_2 \times \cdots \times F_m$, we deduce
\[
D^*\Psi(\bx|A\bx)(w) =
\left\{
\begin{array}{cl}
\{0\} & \mbox{if}~ w_i \in -N_{F_i}(\bx)~ \forall i \\
\emp & \mbox{otherwise}
\end{array}
\right.
\]
and hence
\[
D^*\Phi(\bx|0)(w) =
\left\{
\begin{array}{cl}
-\sum_i w_i & \mbox{if}~ w_i \in -N_{F_i}(\bx)~ \forall i \\
\emp & \mbox{otherwise}.
\end{array}
\right.
\]
From the coderivative formula (\ref{coderivative}) we now obtain
\bmye \label{coderivative2}
\frac{1}{\mbox{reg}\, \Phi(\bx|0)} ~=~
\min \Big\{ \Big\| \sum_i y_i \Big\| : \sum_i \|y_i\|^2 = 1,~ y_i \in N_{F_i}(\bx) \Big\},
\emye
where, following the usual convention, we interpret the right-hand side as $+\infty$ if $N_{F_i}(\bx) = \{0\}$ (or equivalently $\bx \in \inT F_i$) for all $i=1,2,\ldots,m$.  Thus the regularity modulus agrees exactly with the condition modulus that we defined in the previous section:
\[
\mbox{reg}\, \Phi(\bx|0) = \mbox{cond}(F_1,F_2,\ldots,F_m | \bx).
\]
Furthermore, as is well-known \cite{Kru06}, strong regularity is equivalent to the strong metric inequality (\ref{strong}).

\section{Clarke regularity and refinements}
Even more central than strong regularity in variational analysis is the concept of ``Clarke regularity''.  In this section we study a slight refinement, crucial for our development.  In the interest of maintaining as elementary approach as possible, we use the following geometric definition of Clarke regularity.

\begin{defn}[Clarke regularity]
{\rm
A closed set $C \subset \Rn$ is {\em Clarke regular} at a point $\bx \in C$ if, given any $\delta > 0$, any two points $u,z$ near $\bx$ with $z \in C$, and any point $y \in P_C(u)$, satisfy
\[
\ip{z - \bx}{u-y} ~\le~ \delta \|z - \bx\| \cdot \|u-y\|.
\]
}
\end{defn}

\noindent
In other words, the angle between the vectors $z - \bx$ and $u-y$, whenever it is defined, cannot be much less than $\frac{\pi}{2}$ when the points $u$ and $z$ are near $\bx$.

\begin{rem}
{\rm
This property is equivalent to the standard notion of Clarke regularity.  To see this, suppose the property in the definition holds.  Consider any unit vector $v \in N_C(\bx)$, and any unit ``tangent direction'' $w$ to $C$ at $\bx$.  By definition, there exists a sequences $u_r \rt \bx$, $y_r \in P_C(u_r)$, and $z_r \rt \bx$ with $z_r \in C$, such that
\begin{eqnarray*}
\frac{u_r - y_r}{\|u_r - y_r\|} & \rt & v \\ \\
\frac{z_r - \bx}{\|z_r - \bx\|} & \rt & w.
\end{eqnarray*}
By assumption, given any $\e > 0$, for all large $r$ the angle between the two vectors on the left-hand side is at least $\frac{\pi}{2} - \e$, and hence so is the angle between $v$ and $w$.  Thus $\ip{v}{w} \le 0$, so Clarke regularity follows, by \cite[Cor 6.29]{Roc98}.  Conversely, if the property described in the definition fails, then for some $\e > 0$ and some sequences $u_r \rt \bx$, $y_r \in P_C(u_r)$, and $z_r \rt \bx$ with $z_r \in C$, the angle between the unit vectors
\bmye \label{normal-tangent}
\frac{u_r - y_r}{\|u_r - y_r\|} ~~~\mbox{and}~~~ \frac{z_r - \bx}{\|z_r - \bx\|}
\emye
is less than $\frac{\pi}{2} - \e$.  Then any cluster points $v$ and $w$ of the two sequences (\ref{normal-tangent}) are respectively an element of $N_C(\bx)$ and a tangent direction to $C$ at $\bx$, and satisfy $\ip{v}{w} > 0$, contradicting Clarke regularity.
}
\end{rem}

The property we need for our development is an apparently-slight modification of Clarke regularity.

\begin{defn}[super-regularity]
{\rm
A closed set $C \subset \Rn$ is {\em super-regular} at a point $\bx \in C$ if, given any
$\delta > 0$, any two points $u,z$ near $\bx$ with $z \in C$, and any point $y \in P_C(u)$, satisfy
\[
\ip{z - y}{u-y} ~\le~ \delta \|z - y\| \cdot \|u-y\|.
\]
}
\end{defn}

\noindent
In other words, then angle between the vectors $z - y$ and $u-y$, whenever it is defined, cannot be much less than $\frac{\pi}{2}$ when the points $u$ and $z$ are near $\bx$.  

An equivalent statement involves the normal cone.

\begin{prop}[super-regularity and normal angles] \label{super-normal} \hfill
A closed set \mbox{$C \subset \Rn$} is super-regular at a point $\bx \in C$ if and only if, for all
$\delta > 0$, the inequality
\[
\ip{v}{z-y} ~\le~ \delta \|v\| \cdot \|z-y\|
\]
holds for all points $y,z \in C$ near $\bx$ and all normal vectors $v \in N_C(y)$.
\end{prop}

\pf
Super-regularity follows immediately from the normal cone property describe in the proposition, by property (\ref{projections-normals}).  Conversely, suppose the normal cone property fails, so for some $\delta > 0$ and sequences of distinct points $y_r,z_r \in C$ approaching $\bx$ and unit normal vectors
$v_r \in N_C(y_r)$, we have, for all $r=1,2,\ldots$,
\[
\Big\langle v_r , \frac{z_r-y_r}{\|z_r-y_r\|} \Big\rangle ~>~ \delta .
\]

Fix an index $r$.  By definition of the normal cone, there exist sequences of distinct points $u_r^j \rt y_r$ and $y_r^j \in P_C(u_r^j)$ such that
\[
\lim_{j \rt \infty} \frac{u_r^j - y_r^j}{\|u_r^j - y_r^j\|} = v_r.
\]
Since $\lim_j y_r^j = y_r$, we must have, for all large $j$,
\[
\Big\langle \frac{u_r^j - y_r^j}{\|u_r^j - y_r^j\|} ,
\frac{z_r-y_r^j}{\|z_r-y_r^j\|} \Big\rangle ~>~ \delta .
\]
Choose $j$ sufficiently large to ensure both the above inequality and the inequality 
$\|u_r^j - y_r\| < \frac{1}{r}$, and then define points $u'_r = u_r^j$ and $y'_r = y_r^j$.

We now have sequences of points $u'_r,z_r$ approaching $\bx$ with $z_r \in C$, and $y'_r \in P_C(u'_r)$, and satisfying
\[
\Big\langle \frac{u'_r - y'_r}{\|u'_r - y'_r\|} ,
\frac{z_r-y'_r}{\|z_r-y'_r\|} \Big\rangle ~>~ \delta .
\]
Hence $C$ is not super-regular at $\bx$.
\finpf

Super-regularity is a strictly stronger property than Clarke regularity, as the following result and example make clear.

\begin{cor}[super-regularity implies Clarke regularity] \hfill \mbox{}
\mbox{If a closed} set $C \subset \Rn$ is super-regular at a point, then it is also Clarke regular there.
\end{cor}

\pf
Suppose the point in question is $\bx$.  Fix any $\delta > 0$, and set $y = \bx$ in  Proposition \ref{super-normal}.  Then clearly any unit tangent direction $d$ to $C$ at $\bx$ and any unit normal vector $v \in N_C(\bx)$ satisfy 
$\ip{v}{d} \le \delta$.  Since $\delta$ was arbitrary, in fact $\ip{v}{d} \le 0$, so Clarke regularity follows by \cite[Cor 6.29]{Roc98}.
\finpf

\begin{exa}
{\rm
Consider the following function $f \colon \R \rt (-\infty,+\infty]$, taken from an example in \cite{Sha00}:
\[
f(t) = 
\left\{
\begin{array}{ll}
2^r(t-2^r) & (2^r \le t < 2^{r+1},~ r \in {\bf Z}) \\
0          & (t=0) \\
+ \infty   & (t<0).
\end{array}
\right.
\]
The epigraph of this function is Clarke regular at $(0,0)$, but it is not hard to see that it is not super-regular there.  Indeed, a minor refinement of this example (smoothing the set slightly close to the nonsmooth points $(2^r,0)$ and $(2^r,4^{r-1})$) shows that a set can be everywhere Clarke regular, and yet not super-regular.
}
\end{exa}

Super-regularity is a common property:  indeed, it is implied by two well-known properties, that we discuss next.  Following \cite{Roc98}, we say that a set $C \subset \Rn$ is {\em amenable} at a point $\bx \in C$ when 
there exists a neighborhood $U$ of $\bx$, a $C^1$ mapping $G \colon U \rt \R^\ell$, and a closed convex set $D \subset \R^\ell$ containing $G(\bx)$, and satisfying the constraint qualification
\bmye \label{amenable}
N_D(G(\bx)) \cap \ker(\nabla G(\bx)^*) = \{0\},
\emye
such that points $x \in \Rn$ near $\bx$ lie in $C$ exactly when $G(x) \in D$.  In particular, if $C$ is defined by $C^1$ equality and inequality constraints and the Mangasarian-Fromovitz constraint qualification holds at $\bx$, then $C$ is amenable at $\bx$.

\begin{prop}[amenable implies super-regular]
If a closed set $C \subset \Rn$ is amenable at a point in $C$, then it is super-regular there.
\end{prop}

\pf
Suppose the result fails at some point $\bx \in C$.  Assume as in the definition of amenability that, in a neighborhood of $\bx$, the set $C$ is identical with the inverse image $G^{-1}(D)$, where the $C^1$ map $G$ and the closed convex set $D$ satisfy the condition (\ref{amenable}).  Then by definition, for some $\delta > 0$, there are sequences of points $y_r,z_r \in C$ and unit normal vectors $v_r \in N_C(y_r)$ satisfying 
\[
\ip{v_r}{z_r-y_r} > \delta \|z_r - y_r\|,~~~ \mbox{for all $r=1,2,\ldots$}.
\]
It is easy to check the condition 
\[
N_D(G(y_r)) \cap \ker(\nabla G(y_r)^*) = \{0\},
\]
for all large $r$, since otherwise we contradict assumption (\ref{amenable}).  Consequently, using the standard chain rule from \cite{Roc98}, we deduce
\[
N_C(y_r) = \nabla G(y_r)^* N_D(G(y_r)),
\]
so there are normal vectors $u_r \in N_D(G(y_r))$ such that $\nabla G(y_r)^* u_r = v_r$.  The sequence $(u_r)$ must be bounded, since otherwise, by taking a subsequence, we could suppose $\|u_r\| \rt \infty$ and $\|u_r\|^{-1} u_r$ approaches some unit vector $\hat u$, leading to the contradiction
\[
\hat u \in N_D(G(\bx)) \cap \ker(\nabla G(\bx)^*) = \{0\}.
\]

For all large $r$, we now have
\[
\ip{\nabla G(y_r)^* u_r}{z_r-y_r} > \delta \|z_r - y_r\|,
\]
and by convexity we know
\[
\ip{u_r}{G(z_r) - G(y_r)} \le 0.
\]
Adding these two inequalities gives
\[
\ip{u_r}{G(z_r) - G(y_r) - \nabla G(y_r)(z_r-y_r)} < -\delta \|z_r - y_r\|.
\]
But as $r \rt \infty$, the left-hand side is $o(\|z_r - y_r\|)$, since the sequence $(u_r)$ is bounded and $G$ is $C^1$.  This contradiction completes the proof.
\finpf

A rather different refinement of Clarke regularity is the notion of ``prox-regularity''.
Following \cite[Thm 1.3]{Pol00}, we call a set $C \subset \E$ is {\em prox-regular} at a point $\bx \in C$ if the projection mapping $P_C$ is single-valued around $\bx$.  (In this case, clearly $C$ must be locally closed around $\bx$.)  For example, if, in the definition of an amenable set that we gave earlier, we strengthen our assumption on the map $G$ to be $C^2$ rather than just $C^1$, the resulting set must be prox-regular.  On the other hand, the set 
\[
\big\{(s,t) \in \R^2 : t = |s|^{3/2}\big\}
\]
is amenable at the point $(0,0)$ (and hence super-regular there), but is not prox-regular there.

\begin{prop}[prox-regular implies super-regular] \label{prox-super}
If a closed set $C \subset \Rn$ is prox-regular at a point in $C$, then it is super-regular there.
\end{prop}

\pf
If the results fails at $\bx \in C$, then for some constant $\delta > 0$, there exist sequences of points $y_r,z_r \in C$ converging to the point $\bx$, and a sequence of normal vectors $v_r \in N_C(y_r)$ satisfying the inequality
\[
\ip{v_r}{z_r-y_r} > \delta\|v_r\| \cdot \|z_r-y_r\|.
\]
By \cite[Proposition 1.2]{Pol00}, there exist constants $\epsilon, \rho > 0$ such that
\[
\Big\langle \frac{\epsilon}{2\|v_r\|} v_r , z_r - y_r \Big\rangle
\le \frac{\rho}{2} \|z_r - y_r\|^2
\]
for all large $r$.  This gives a contradiction, since
$\|z_r - y_r\| \le \frac{\delta \epsilon}{\rho}$ eventually.
\finpf

Super-regularity is related to various other notions in the literature.  We end this section with a brief digression to discuss these relationships.  First note the following equivalent definition, which is an immediate consequence of Proposition \ref{super-normal}, and which gives an alternate proof of Proposition \ref{prox-super} via ``hypomonotonicity'' of the truncated normal cone mapping $x \mapsto N_C(x) \cap B$ for prox-regular sets $C$ \cite[Thm 1.3]{Pol00}.

\begin{cor}[approximate monotonicity]
\hfill A closed set \mbox{$C \subset \Rn$} is super-regular at a point $\bx \in C$ if and only if, for all $\delta > 0$, the inequality
\[
\ip{v-w}{y-z} ~\ge~ -\delta \|y-z\|
\]
holds for all points $y,z \in C$ near $\bx$ and all normal vectors $v \in N_C(y) \cap B$ and $w \in N_C(z) \cap B$.
\end{cor}

\noindent
If we replace the normal cone $N_C$ in the property described in the result above by its convex hull, the ``Clarke normal cone'', we obtain a stronger property, called ``subsmoothness'' in \cite{Aus04}.  Similar proofs to those above show that, like super-regularity, subsmoothness is a consequence of either amenability or prox-regularity.  However, submoothness is strictly stronger than super-regularity.  To see this, consider the graph of the function $f \colon \R \rt \R$ defined by the following properties: $f(0) = 0$, $f(2^r) = 4^r$ for all integers $r$, $f$ is linear on each interval $[2^r,2^{r+1}]$, and $f(t) = f(-t)$ for all 
$t \in \R$.  The graph of $f$ is super-regular at $(0,0)$, but is not subsmooth there.

In a certain sense, however, the distinction between subsmoothness and super-regularity is slight.  Suppose the set $F$ is super-regular at every point in $F \cap U$, for some open set 
$U \subset \Rn$.  Since super-regularity implies Clarke regularity, the normal cone and Clarke normal cone coincide throughout $F \cap U$, and hence $F$ is also subsmooth throughout 
$F \cap U$.  In other words, ``local'' super regularity coincides with ``local'' subsmoothness, which in turn, by \cite[Thm 3.16]{Aus04} coincides with the ``first order Shapiro property'' \cite{Sha94} (also called ``near convexity'' in \cite{Sha93}) holding locally.

\section{Alternating projections with nonconvexity}

Having reviewed or developed over the last few sections the key variational-analytic properties that we need, we now turn to projection algorithms.
In this section we develop our convergence analysis of the method of alternating projections.
The following result is our basic tool, guaranteeing conditions under
which the method of alternating projections converges linearly.
For flexibility, we state it in a rather technical manner.
For clarity, we point out afterward that the two main conditions,
(\ref{condition}) and (\ref{condition'}), are guaranteed in applications
via assumptions of strong regularity and super-regularity (or in particular, amenability or prox-regularity) respectively.

\begin{thm}[linear convergence of alternating projections] \label{main} \hfill \mbox{}
Consider the closed sets $F,C \subset \E$, and a point $\bx \in F$.  Fix any constant $\epsilon > 0$.  Suppose for some constant $c' \in (0,1)$, the following condition holds:
\bmye \label{condition}
\left.
\begin{array}{cc}
x \in F \cap (\bx + \epsilon B), &  u \in -N_F(x) \cap B  \\
y \in C \cap (\bx + \epsilon B), &  v \in N_C(y) \cap B
\end{array}
\right\}
~~\Rightarrow~~ \ip{u}{v} \le c'.
\emye
Suppose furthermore for some constant $\delta \in [0,\frac{1-c'}{2})$ the following condition holds:
\bmye \label{condition'}
\left.
\begin{array}{rcl}
y,z & \in & C \cap (\bx + \epsilon B) \\
  v & \in & N_C(y) \cap B
\end{array}
\right\}
~~\Rightarrow~~ \ip{v}{z-y} \le \delta\|z-y\|.
\emye
Define a constant $c = c'+2\delta < 1$.  Then for any initial point $x_0 \in C$ satisfying
$\|x_0 - \bx\| \le \frac{1-c}{4} \epsilon$, any sequence of alternating projections on the sets $F$ and $C$,
\[
x_{2n+1} \in P_F(x_{2n}) ~~ \mbox{and} ~~ x_{2n+2} \in P_C(x_{2n+1}) ~~~ (n=0,1,2,\ldots)
\]
must converge with R-linear rate $\sqrt{c}$ to a point
$\hat x \in F \cap C$ satisfying the inequality $\| \hat x - x_0 \| \le \frac{1+c}{1-c}\|x_0 - \bx\|$.
\end{thm}

\pf
First note, by the definition of the projections we have
\bmye \label{easy2'}
\|x_{2n+3} - x_{2n+2}\| ~\le~ \|x_{2n+2} - x_{2n+1}\| ~\le~ \|x_{2n+1} - x_{2n}\|.
\emye
Clearly we therefore have
\bmye \label{easy4}
\|x_{2n+2} - x_{2n}\| \le 2\|x_{2n+1} - x_{2n}\|.
\emye
We next claim
\bmye \label{implication'}
\left.
\begin{array}{c}
\|x_{2n+1} - \bx\| \le \frac{\epsilon}{2} ~~\mbox{and} \\
\|x_{2n+1} - x_{2n}\| \le \frac{\epsilon}{2}
\end{array}
\right\}
~~ \Rightarrow ~~
\|x_{2n+2} - x_{2n+1}\| \le c \|x_{2n+1} - x_{2n}\|.
\emye
To see this, note that if $x_{2n+2} = x_{2n+1}$, the result is trivial, and if $x_{2n+1} = x_{2n}$ then $x_{2n+2} = x_{2n+1}$ so again the result is trivial.  Otherwise, we have
\[
\frac{x_{2n} - x_{2n+1}}{\|x_{2n} - x_{2n+1}\|} \in N_F(x_{2n+1}) \cap B
\]
while
\[
\frac{x_{2n+2} - x_{2n+1}}{\|x_{2n+2} - x_{2n+1}\|} \in -N_C(x_{2n+2}) \cap B.
\]
Furthermore, using inequality (\ref{easy2'}), the left-hand side of the implication (\ref{implication'}) ensures
\begin{eqnarray*}
\|x_{2n+2} - \bx\|
& \le &
\|x_{2n+2} - x_{2n+1}\| + \|x_{2n+1} - \bx\| \\
& \le &
\|x_{2n+1} - x_{2n}\| + \|x_{2n+1} - \bx\| ~\le~ \epsilon.
\end{eqnarray*}
Hence, by assumption (\ref{condition}) we deduce
\[
\Big\langle
\frac{x_{2n} - x_{2n+1}}{\|x_{2n} - x_{2n+1}\|} \:,\:
\frac{x_{2n+2} - x_{2n+1}}{\|x_{2n+2} - x_{2n+1}\|}
\Big\rangle
~\le~ c',
\]
so
\[
\ip{x_{2n} - x_{2n+1}}{x_{2n+2} - x_{2n+1}} ~\le~
c' \|x_{2n} - x_{2n+1}\| \cdot \|x_{2n+2} - x_{2n+1}\|.
\]
On the other hand, by assumption (\ref{condition'}) we know
\begin{eqnarray*}
\ip{x_{2n} - x_{2n+2}}{x_{2n+1} - x_{2n+2}}
& \le &
\delta \|x_{2n} - x_{2n+2}\| \cdot \|x_{2n+1} - x_{2n+2}\|  \\
& \le &
2\delta \|x_{2n} - x_{2n+1}\| \cdot \|x_{2n+2} - x_{2n+1}\|,
\end{eqnarray*}
using inequality (\ref{easy4}).
Adding this inequality to the previous inequality then gives the right-hand side of (\ref{implication'}), as desired.

Now let $\alpha = \|x_0 - \bx\|$.
We will show by induction the inequalities
\begin{myeqnarray}
\|x_{2n+1} - \bx\|       & \le & 2\alpha \frac{1-c^{n+1}}{1-c} ~<~ \frac{\epsilon}{2}
                                                                    \label{claim1'} \\
\|x_{2n+1} - x_{2n} \|   & \le & \alpha c^n ~<~ \frac{\epsilon}{2}  \label{claim2'}  \\
\|x_{2n+2} - x_{2n+1} \| & \le & \alpha c^{n+1}.                    \label{claim3'}
\end{myeqnarray}

Consider first the case $n=0$.  Since $x_1 \in P_F(x_0)$ and $\bx \in F$, we deduce
$\|x_1 - x_0\| \le \|\bx - x_0\| = \alpha < \epsilon/2$, which is inequality (\ref{claim2'}).  Furthermore,
\[
\|x_1 - \bx\| \le \|x_1 - x_0\| + \|x_0 - \bx\|  \le 2\alpha < \frac{\epsilon}{2},
\]
which shows inequality (\ref{claim1'}).  Finally, since $\|x_1 - x_0\| < \epsilon/2$ and
$\|x_1 - \bx\| < \epsilon/2$, the implication (\ref{implication'}) shows
\[
\| x_2 - x_1 \| \le c\|x_1 - x_0\| \le c\|\bx - x_0\| = c\alpha,
\]
which is inequality (\ref{claim3'}).

For the induction step, suppose inequalities (\ref{claim1'}), (\ref{claim2'}), and (\ref{claim3'}) all hold for some $n$.  Inequalities (\ref{easy2'}) and (\ref{claim3'}) imply
\bmye \label{nclaim2'}
\|x_{2n+3} - x_{2n+2}\| \le \alpha c^{n+1} < \frac{\epsilon}{2}.
\emye
We also have, using inequalities (\ref{nclaim2'}), (\ref{claim3'}), and (\ref{claim1'})
\begin{eqnarray*}
\|x_{2n+3} - \bx\|
& \le &
\|x_{2n+3} - x_{2n+2}\| + \|x_{2n+2} - x_{2n+1}\| + \|x_{2n+1} - \bx\| \\
& \le &
\alpha c^{n+1} + \alpha c^{n+1} + 2\alpha \frac{1-c^{n+1}}{1-c},
\end{eqnarray*}
so
\bmye \label{nclaim1'}
\|x_{2n+3} - \bx\| \le 2\alpha \frac{1-c^{n+2}}{1-c} < \frac{\epsilon}{2}.
\emye
Now implication (\ref{implication'}) with $n$ replaced by $n+1$ implies
\[
\|x_{2n+4} - x_{2n+3}\| \le c \|x_{2n+3} - x_{2n+2}\|,
\]
and using inequality (\ref{nclaim2'}) we deduce
\bmye \label{nclaim3'}
\|x_{2n+4} - x_{2n+3}\| \le \alpha c^{n+2}.
\emye
Since inequalities (\ref{nclaim1'}), (\ref{nclaim2'}), and (\ref{nclaim3'}) are exactly inequalities (\ref{claim1'}), (\ref{claim2'}), and (\ref{claim3'}) with $n$ replaced by $n+1$, the induction step is complete and our claim follows.

We can now easily check that the sequence $(x_k)$ is Cauchy and therefore converges.  To see this, note for any integer $n=0,1,2,\ldots$ and any integer $k > 2n$, we have
\begin{eqnarray*}
\|x_k - x_{2n}\|
& \le &
\sum_{j=2n}^{k-1} \|x_{j+1} - x_j\| \\
& \le &
\alpha( c^n + c^{n+1} + c^{n+1} + c^{n+2} + c^{n+2} + \cdots )
\end{eqnarray*}
so
\[
\|x_k - x_{2n}\| \le \alpha c^n \frac{1+c}{1-c},
\]
and a similar argument shows
\bmye \label{e:J1}
\|x_{k+1} - x_{2n+1}\| \le \frac{2\alpha c^{n+1}}{1-c}.
\emye
Hence $x_k$ converges to some point $\hat x \in \E$, and for all $n=0,1,2,\ldots$ we have
\bmye \label{e:J2}
\|\hat x - x_{2n}\| \le \alpha c^n \frac{1+c}{1-c}
~~~\mbox{and}~~~
\|\hat x - x_{2n+1}\| \le \frac{2\alpha c^{n+1}}{1-c}.
\emye
We deduce that the limit $\hat x$ lies in the intersection $F \cap C$ and satisfies the inequality
$\|\hat x - x_0\| \le \alpha \frac{1+c}{1-c}$, and furthermore that the convergence is R-linear with rate $\sqrt{c}$, which completes the proof.
\finpf

To apply Theorem \ref{main} to alternating projections between
a closed and a super-regular set, we make use of the key geometric
property of super-regular sets (Proposition \ref{super-normal}):  at any point near a point where a set is super-regular the angle between any normal vector and the direction
to any nearby point in the set cannot be much less than $\frac{\pi}{2}$.

We can now prove our key result.

\begin{thm}[alternating projections with a super-regular set] \label{key} \hfill \mbox{}
Consider closed sets $F,C \subset \E$ and a point $\bx \in F \cap C$.  Suppose $C$ is super-regular at $\bx$ (as holds, for example, if it is amenable or prox-regular there).  Suppose furthermore that $F$ and $C$ have strongly regular intersection at $\bx$:  that is, the condition
\[
N_F(\bx) \cap -N_C(\bx) = \{0\}
\]
holds, or equivalently, the constant
\bmye \label{barc}
\bar c ~=~ \max \Big\{ \ip{u}{v} : u \in N_F(\bx) \cap B,~ v \in -N_C(\bx) \cap B \Big\}
\emye
is strictly less than one. Fix any constant $c \in (\bar c,1)$.  Then, for any initial point
$x_0 \in C$ close to $\bx$, any sequence of iterated projections
\[
x_{2n+1} \in P_F(x_{2n}) ~~ \mbox{and} ~~ x_{2n+2} \in P_C(x_{2n+1}) ~~~ (n=0,1,2,\ldots)
\]
must converge to a point in $F \cap C$ with R-linear rate $\sqrt{c}$.
\end{thm}

\pf
 Let us show first the equivalence between $\bar c <1$ and strong regularity.
The compactness of the intersections between normal cones and the unit ball
guarantees the existence of $u$ and $v$ achieving the maximum in \eqref{barc}.
Observe then that
\[
\langle{u},{v}\rangle\leq \norm{u}\,\norm{v}\leq 1.
\]
The cases of equality in the Cauchy-Schwarz inequality permits to write
\[
\bar c =1  \iff \mbox{$u$ and $v$ are colinear} \iff N_F(\bx) \cap -N_C(\bx) \neq \{0\},
\]
which corresponds to the desired equivalence.

Fix now any constant $c' \in (\bar c, c)$ and define $\delta = \frac{c-c'}{2}$.
To apply Theorem \ref{main}, we just need to check the existence of a constant $\e > 0$ such that conditions (\ref{condition}) and (\ref{condition'}) hold.  Condition (\ref{condition'}) holds for all small $\e > 0$, by Proposition \ref{super-normal}.  On the other hand, if condition (\ref{condition}) fails for all small $\e > 0$, then there exist sequences of points $x_r \rt \bx$ in the set $F$ and $y_r \rt \bx$ in the set $C$, and sequences of vectors
$u_r \in -N_F(x_r) \cap B$ and $v_r \in N_C(y_r) \cap B$, satisfying $\ip{u_r}{v_r} > c'$.  After taking a subsequences, we can suppose $u_r$ approaches some vector
$u \in -N_F(\bx) \cap B$ and $v_r$ approaches some vector $v \in N_C(\bx) \cap B$, and then $\ip{u}{v} \ge c' > \bar c$, contradicting the definition of the constant $\bar c$.
\finpf

\begin{cor}[improved convergence rate] \label{improved}
With the assumptions of Theorem \ref{key}, suppose the set $F$ is also super-regular at $\bx$.  Then the alternating projection sequence converges with R-linear rate $c$.
\end{cor}

\pf
Inequality (\ref{implication'}), and its analog when the roles of $F$ and $C$ are interchanged, together show
\[
\| x_{k+1} - x_k \| \le c \|x_k - x_{k-1}\|
\]
for all large $k$, and the result then follows easily, using an argument analogous to that at the end of the proof of Theorem \ref{main}.
\finpf

In the light of our discussion in the previous section, the strong regularity assumption of Theorem \ref{key} is equivalent to the metric regularity at $\bx$ for $0$ of the set-valued mapping $\Psi \colon \E \tto \E^2$ defined by
\[
\Psi(x) = (F-x) \times (C-x),~~ \mbox{for}~ x \in \E.
\]
Using equation (\ref{coderivative2}), the regularity modulus is determined by
\[
\frac{1}{\mbox{\rm reg}\,\Psi(\bx|0)} ~=~
\min \Big\{ \|u+v\| : u \in N_F(\bx),~ v \in N_C(\bx),~
\|u\|^2 + \|v\|^2 = 1 \Big\},
\]
and a short calculation then shows
\bmye \label{reg-c}
\mbox{\rm reg}\,\Psi(\bx|0) = \frac{1}{\sqrt{1- \bar c}}.
\emye
The closer the constant $\bar c$ is to one, the larger the regularity modulus.  We have shown that $\bar c$ also controls the speed of linear convergence for the method of alternating projections applied to the sets $F$ and $C$.

Inevitably, Theorem \ref{key} concerns {\em local} convergence:  it relies on finding an initial point $x_0$ sufficiently close to a point of strongly regular intersection.  How might we find such a point?

One natural context in which to pose this question is that of sensitivity analysis.  Suppose we already know a point of strongly regular intersection of a closed set and a subspace, but now want to find a point in the intersection of two slight perturbations of these sets.  The following result shows that, starting from the original point of intersection, the method of alternating projections will converge linearly to the new intersection.

\begin{thm}[perturbed intersection] \label{perturbed}
With the assumptions of Theorem \ref{key}, for any small vector $d \in \E$, the method of alternating projections applied to the sets $d+F$ and $C$, with the initial point
$\bx \in C$, will converge with R-linear rate $\sqrt{c}$ to a point $\tilde x \in (d+F) \cap C$ satisfying $\|\tilde x - \bar x\| \le \frac{1+c}{1-c} \|d\|$.
\end{thm}

\pf
As in the proof of Theorem \ref{key}, if we fix any constant $c' \in (\bar c, c)$ and define
$\delta  = \frac{c-c'}{2}$, then there exists a constant $\e > 0$ such that conditions (\ref{condition}) and (\ref{condition'}) hold.  Suppose the vector $d$ satisfies
\[
\|d\| \le \frac{(1-c)\e}{8} < \frac{\e}{2}.
\]
Since
\begin{eqnarray*}
\lefteqn{
y \in (C-d) \cap (\bx + \frac{\e}{2} B) ~\mbox{and}~ v \in N_{C-d}(y)
}\\
& &
~~\Rightarrow~~ y+d \in C \cap (\bx + \e B)  ~\mbox{and}~ v \in N_C(y+d),
\end{eqnarray*}
we deduce from condition (\ref{condition}) the implication
\[
\left.
\begin{array}{cc}
x \in F \cap (\bx + \frac{\epsilon}{2} B), &  u \in -N_F(x) \cap B  \\
y \in (C-d) \cap (\bx + \frac{\epsilon}{2} B), &  v \in N_{C-d}(y) \cap B
\end{array}
\right\}
~~\Rightarrow~~ \ip{u}{v} \le c'.
\]
Furthermore, using condition (\ref{condition'}) we deduce the implication
\begin{eqnarray*}
\lefteqn{
y,z \in (C-d) \cap (\bx + \frac{\e}{2} B) ~\mbox{and}~ v \in N_{C-d}(y) \cap B
}\\
& &
~~\Rightarrow~~ y+d, z+d \in C \cap (\bx + \e B)  ~\mbox{and}~ v \in N_C(y+d) \cap B, \\
& &
~~\Rightarrow~~ \ip{v}{z-y} \le \delta \|z-y\|.
\end{eqnarray*}
We can now apply Theorem \ref{main} with the set $C$ replaced by $C-d$ and the constant $\e$ replaced by $\frac{\e}{2}$.  We deduce that the method of alternating projections on the sets $F$ and $C-d$, starting at the point $x_0 = \bx - d \in C-d$, converges with R-linear rate $\sqrt{c}$ to a point $\hat x \in F \cap (C-d)$ satisfying the inequality
$\| \hat x - x_0 \| \le \frac{1+c}{1-c} \|x_0 - \bx\|$.  The theorem statement then follows by translation.
\finpf

Lack of convexity notwithstanding, more structure sometimes implies that the method of alternating projections converges Q-linearly, rather than just R-linearly, on a neighborhood of point of strongly regular intersection of two closed sets.  One example is the case of two manifolds \cite{Lew06Alternating}.

\section{Inexact alternating projections}
Our basic tool, the method of alternating projections for a super-regular set $C$ and an arbitrary closed set $F$, is a conceptual algorithm that may be challenging to realize in practice.  We might reasonably consider the case of exact projections on the super-regular set $C$:  for example, in the next section, for the method of averaged projections, $C$ is a subspace and computing projections is trivial.  However, projecting onto the set $F$ may be much harder, so a more realistic analysis allows relaxed projections.

We sketch one approach.  Given two iterates $x_{2n-1} \in F$ and $x_{2n} \in C$, a necessary condition for the new iterate $x_{2n+1}$ to be an exact projection on $F$, that is
$x_{2n+1} \in P_F(x_{2n})$, is
\[
\|x_{2n+1} - x_{2n}\| \le \|x_{2n} - x_{2n-1}\|
~~\mbox{and}~~
x_{2n} - x_{2n+1} \in N_F(x_{2n+1}).
\]
In the following result we assume only that we choose the iterate $x_{2n+1}$ to satisfy a relaxed version of this condition, where we replace the second part by the assumption that the distance
\[
d_{N_F(x_{2n+1})} \Big( \frac{x_{2n} - x_{2n+1}}{\|x_{2n} - x_{2n+1}\|} \Big)
\]
from the normal cone at the iterate to the normalized direction of the last step is {\em small}.

\begin{thm}[inexact alternating projections]
\label{inexact}
With the assumptions of Theorem \ref{key}, fix any constant $\epsilon < \sqrt{1-c^2}$, and consider the following inexact alternating projection iteration.  Given any initial points $x_0 \in C$ and $x_1 \in F$, for $n=1,2,3,\ldots$ suppose
\[
x_{2n} \in P_C(x_{2n-1})
\]
and $x_{2n+1} \in F$ satisfies
\[
\|x_{2n+1} - x_{2n}\| \le \|x_{2n}-x_{2n-1}\|
~~\mbox{and}~~
d_{N_F(x_{2n+1})} \Big( \frac{x_{2n} - x_{2n+1}}{\|x_{2n} - x_{2n+1}\|} \Big) \le \epsilon.
\]
Then, providing $x_0$ and $x_1$ are close to $\bx$, the iterates converge to a point in
$F \cap C$ with R-linear rate
\[
\sqrt{ c\sqrt{1-\epsilon^2} + \epsilon \sqrt{1-c^2} } ~<~ 1.
\]
\end{thm}

\linespace
\noindent
{\bf Sketch proof.}
Once again as in the proof of Theorem \ref{key}, we fix any constant $c' \in (\bar c, c)$ and define $\delta  = \frac{c-c'}{2}$, so there exists a constant $\e > 0$ such that conditions (\ref{condition}) and (\ref{condition'}) hold.
Define a vector
\[
z = \frac{x_{2n} - x_{2n+1}}{\|x_{2n} - x_{2n+1}\|}.
\]
By assumption, there exists a vector $w \in N_F(x_{2n+1})$ satisfying $\|w-z\| \le \epsilon$.  Some elementary manipulation then shows that the unit vector $\hat w = \|w\|^{-1} w$ satisfies
\[
\ip{\hat w}{z} \ge \sqrt{1-\epsilon^2}.
\]
As in the proof of Theorem \ref{main}, assuming inductively that $x_{2n+1}$ is close to both $\bx$ and $x_{2n} $, since $\hat w \in N_F(x_{2n+1})$, and
\[
u = \frac{x_{2n+2}-x_{2n+1}}{\|x_{2n+2}-x_{2n+1}\|} \in -N_C(x_{2n+2}) \cap B,
\]
we deduce
\[
\ip{\hat w}{u} \le c'.
\]

We now see that, on the unit sphere, the arc distance between the unit vectors $\hat w$ and $z$ is no more than $\arccos(\sqrt{1-\epsilon^2})$, whereas the arc distance between $\hat w$ and the unit vector $u$ is at least $\arccos c'$.  Hence by the triangle inequality, the arc distance between $z$ and $u$ is at least
\[
\arccos c' - \arccos(\sqrt{1-\epsilon^2}),
\]
so
\[
\ip{z}{u} ~\le~ \cos \Big( \arccos c' - \arccos(\sqrt{1-\epsilon^2}) \Big)
~=~
c'\sqrt{1-\epsilon^2} + \epsilon \sqrt{1-c'^2}.
\]
Some elementary calculus shows that the quantity on the right-hand side is strictly less than one.  Again as in the proof of Theorem \ref{main}, this inequality shows, providing $x_0$ is close to $\bx$, the inequality
\[
\|x_{2n+2} - x_{2n+1}\| ~\le~
\Big( c\sqrt{1-\epsilon^2} + \epsilon \sqrt{1-c^2} \Big) \|x_{2n+1} - x_{2n}\|,
\]
and in conjunction with the inequality
\[
\|x_{2n+1} - x_{2n}\| ~\le~
\|x_{2n} - x_{2n-1}\|,
\]
this suffices to complete the proof by induction.
\finpf

\section{Local convergence for averaged projections}

 We now return to the problem of finding a point in the intersection of
several closed sets using the method of averaged projections.  The results of the
previous section are applied to the method of averaged projections via the well-known
reformulation of the algorithm as alternating projections on a product space.
This leads to the main result of this section, Theorem \ref{lin-con}, which shows linear convergence in a neighborhood of any point of strongly regular intersection,
at a rate governed by the associated regularity modulus.

We begin with a characterization of strongly regular intersection, relating the condition modulus with a generalized notion of angle for several sets.  Such notions, for collections of convex sets, have also been studied recently in the context of projection algorithms in \cite{Deu06a,Deu06b}.

\begin{prop}[variational characterization \hspace{-.8pt}of\hspace{-.8pt} strong regularity\!] \label{several}
\hfill \mbox{}
Closed sets $F_1,F_2,\ldots,F_m \subset \E$ have strongly regular intersection at a point $\bx \in \cap_i F_i$ if and only if the optimal value $\bar c$ of the optimization problem
\begin{eqnarray*}
\mbox{maximize} & & \sum_i \ip{u_i}{v_i} \\
\mbox{subject to} & & \sum_i \|u_i\|^2 \le 1 \\
                  & & \sum_i \|v_i\|^2 \le 1 \\
                  & & \sum_i u_i = 0 \\
                  & & u_i \in \E,~ v_i \in N_{F_i}(\bx) ~~~(i=1,2,\ldots,m)
\end{eqnarray*}
is strictly less than one. Indeed, we have
\bmye \label{cbar}
{\bar c}^2 ~=~
\left\{
\begin{array}{cl}
0 & (\bx \in \cap_i \inT F_i) \\
1 - \displaystyle{\frac{1}{m \cdot \mbox{\rm cond}^2(F_1,F_2,\ldots,F_m | \bx)}} & (\mbox{otherwise}).
\end{array}
\right.
\emye
\end{prop}

\pf
When $\bx \in \cap_i \inT F_i$, the result follows by definition.  Henceforth, we therefore rule out that case.

For any vectors $u_i,v_i \in \E$ ($i=1,2,\ldots,m$), by Lagrangian duality and differentiation we obtain
\begin{eqnarray*}
\lefteqn{
\max_{u_i} \Big\{ \sum_i \ip{u_i}{v_i} : \sum_i \|u_i\|^2 \le 1,~ \sum_i u_i = 0 \Big\}
} \\
& & =
\min_{\lambda \in \R_+,~ z \in \E} \max_{u_i}
\Big\{
\sum_i \ip{u_i}{v_i} + \frac{\lambda}{2} \Big( 1 - \sum_i \|u_i\|^2 \Big)
+ \ip{z}{\sum_i u_i} \Big\} \\
& & =
\min_{\lambda \in \R_+,~ z \in \E}
\Big\{
\frac{\lambda}{2} + \sum_i \max_{u_i} \Big\{ \ip{u_i}{v_i + z} - \frac{\lambda}{2} \|u_i\|^2 \Big\}
\Big\} \\
& & =
\min_{\lambda > 0,~ z \in \E}
\Big\{
\frac{\lambda}{2} + \frac{1}{2\lambda} \sum_i \|v_i + z\|^2
\Big\} \\
& & =
\min_{z \in \E} \sqrt{\sum_i \|v_i + z\|^2} \\
& & = \sqrt{\sum_{i=1}^m \Big\|v_i - \frac{1}{m} \sum_j v_j \Big\|^2}
~=~ \sqrt{ \sum_i \|v_i\|^2 - \frac{1}{m}\Big\|\sum_i v_i \Big\|^2 }.
\end{eqnarray*}
Consequently, ${\bar c}^2$ is the optimal value of the optimization problem
\begin{eqnarray*}
\mbox{maximize}   & & \sum_i \|v_i\|^2 - \frac{1}{m}\Big\|\sum_i v_i \Big\|^2 \\
\mbox{subject to} & & \sum_i \|v_i\|^2 \le 1 \\
                  & & v_i \in N_{F_i}(\bx) ~~~(i=1,2,\ldots,m).
\end{eqnarray*}
By homogeneity, the optimal solution must occur when the inequality constraint is active, so we obtain an equivalent problem by replacing that constraint by the corresponding equation.   By \ref{coderivative2} an the definition of the condition modulus
it follows that the optimal value of this new problem is
\[
1 - \frac{1}{m \cdot \mbox{\rm cond}^2(F_1,F_2,\ldots,F_m | \bx)}
\]
as required.
\finpf

\begin{thm}[linear convergence of averaged projections] \label{lin-con} \hfill
Suppose closed sets $F_1,F_2,\ldots,F_m \subset \E$ have strongly regular intersection at a point $\bx \in \cap_i F_i$.  Define a constant $\bar c \in [0,1)$ by equation (\ref{cbar}), and fix any constant $c \in (\bar c ,1)$.  Then for any initial point $x_0 \in \E$ close to $\bx$, the method of averaged projections converges to a point in the intersection $\cap_i F_i$, with R-linear rate $c$ (and if each set $F_i$ is super-regular at $\bx$, or in particular, prox-regular or amenable there, then the convergence rate is $c^2$).  Furthermore, for any small perturbations $d_i \in \E$ for $i=1,2,\ldots,m$, the method of averaged projections applied to the sets $d_i + F_i$, with the initial point $\bx$, converges linearly to a nearby point in the intersection, with R-linear rate~$c$.
\end{thm}

\pf
In the product space $\E^m$ with the inner product
\[
\ip{(u_1,u_2,\ldots,u_m)}{(v_1,v_2,\ldots,v_m)} ~=~ \sum_i \ip{u_i}{v_i},
\]
we consider the closed set
\[
F = \prod_i F_i
\]
and the subspace
\[
L =\{ Ax : x \in \E \},
\]
where the linear map $A \colon \E \rt \E^m$ is defined by $Ax = (x,x,\ldots,x)$.
Notice $A\bx \in F \cap L$,
and it is easy to check
\[
N_F(A\bx) = \prod_i N_{F_i}(\bx)
\]
and
\[
L^{\perp} ~=~ \Big\{ (u_1,u_2,\ldots,u_m) : \sum_i u_i = 0 \Big\}.
\]
Hence $F_1,F_2,\ldots,F_m$ have strongly regular intersection at $\bx$ if and only if $F$ and $L$ have strongly regular intersection at the point $A\bx$.  This latter property is equivalent to the constant $\bar c$ defined in Theorem \ref{key} (with $C=L$) being strictly less than one.  But that constant agrees exactly with that defined by equation (\ref{cbar}), so
we show next that we can apply Theorem \ref{key} and Theorem \ref{perturbed}.

To see this note that, for any point $x \in \E$, we have the equivalence
\[
(z_1,z_2,\ldots,z_m) \in P_F(Ax) ~~\Leftrightarrow~~
z_i \in P_{F_i}(x) ~~(i=1,2,\ldots,m).
\]
Furthermore a quick calculation shows, for any $z_1,z_2,\ldots,z_m \in \E$,
\[
P_L(z_1,z_2,\ldots,z_m) ~=~ \frac{1}{m} (z_1 + z_2 + \cdots + z_m).
\]
Hence in fact the method of averaged projections for the sets $F_1,F_2,\ldots,F_m$, starting at an initial point $x_0$, is essentially identical with the method of alternating projections for the sets $F$ and $L$, starting at the initial point $Ax_0$.  If  $x_0,x_1,x_2,\ldots$ is a possible sequence of iterates for the former method, then a possible sequence of even iterates for the latter method is $Ax_0,Ax_1,Ax_2,\ldots$.
For $x_0$ close to $\bx$, this latter sequence must converge to a point $A\hat x \in F \cap L$  with R-linear rate $c$, by Theorem \ref{key} and its corollary.
Thus the sequence $x_0,x_1,x_2,\ldots$ converges
to $\hat x \in \cap_i F_i$ at the same linear rate.  When each of the sets $F_i$ is super-regular at $\bx$, it is easy to check that the Cartesian product $F$ is super-regular at $A\bx$, so the rate is $c^2$.  The last part of the theorem follows from Theorem \ref{perturbed}.
\finpf

\noindent
Applying Theorem \ref{inexact} to the
product-space formulation
of averaged projections shows in a similar fashion that an inexact variant of
 the method of averaged projections will also converge linearly.

\begin{rem}[strong regularity and local extremality]
{\rm
  We notice, in the language of \cite{Mor06}, that
  we have proved algorithmically
  that if closed sets have strongly regular intersection at a point,
  then that point is not ``locally extremal''.
}
\end{rem}

\begin{rem}[alternating versus averaged projections]\label{alt-vs-av}
\hfill
{\rm
  \mbox{Consider a} feasibility problem involving two super-regular sets $F_1$ and $F_2$.
  Assume that strong regularity holds at $\bar x \in F_1\cap F_2$
  and set $\kappa = \mbox{{\rm cond}}(F_1,F_2|\bar x)$.
  Theorem \ref{lin-con} gives a bound on the rate of convergence
of the method of averaged projections as
 \[
  r_{\mbox{{\footnotesize{\rm av}}}} \leq 1-\frac{1}{2\kappa^2}\ .
  \]
Notice that each iteration involves two projections:  one onto each of the sets $F_1$ and $F_2$.  On the other hand, Corollary \ref{improved} and (\ref{reg-c}) give a bound on the
rate of convergence of the method of alternating projections as
  \[
  r_{\mbox{{\footnotesize{\rm alt}}}} \leq 1-\frac{1}{\kappa^2}\ ,
  \]
and each iteration involves just one projection.
  Thus we note that our bound on the rate of alternating projections $r_{\mbox{{\footnotesize{\rm alt}}}}$
  is always better than the bound on the rate of averaged projection $r_{\mbox{{\footnotesize{\rm av}}}}$.
At least from the perspective of this analysis, averaged projections seems to have no advantage over alternating projections,
although our proof of linear convergence for alternating projections needs a super-regularity assumption not necessary in the case of
averaged projections.
}
\end{rem}

\section{Prox-regularity and averaged projections}\label{s:dist}
If we assume that the sets $F_1,F_2,\ldots,F_m$ are prox-regular, then we can refine our understanding of local convergence for the method of averaged projections using a completely different approach, explored in this section.

\begin{prop}
Around any point $\bx$ at which the set $F \subset \E$ is prox-regular, the squared distance to $F$ is continuously differentiable, and its gradient $\nabla d_F^2 = 2(I-P_F)$ has Lipschitz constant $2$.
\end{prop}

\pf
 This result corresponds essentially to \cite[Prop 3.1]{Pol00},
which yields the smoothness of $d_F^2$ together with the gradient formula.
This proof of this proposition also shows that for any small $\delta > 0$,
all points $x_1,x_2 \in \E$ near $\bx$ satisfy the inequality
\[
\ip{x_1-x_2}{P_F(x_1) - P_F(x_2)} \ge (1 - \delta)\|P_F(x_1) - P_F(x_2)\|^2
\]
(see ``Claim'' in \cite[p.~5239]{Pol00}).  Consequently we have
\begin{eqnarray*}
\lefteqn{
\|(I-P_F)(x_1) - (I-P_F)(x_2)\|^2 - \|x_1 - x_2\|^2
} \\
& & = \|(x_1 - x_2) - (P_F(x_1) - P_F(x_2))\|^2 - \|x_1 - x_2\|^2 \\
& & = -2\ip{x_1 - x_2}{P_F(x_1) - P_F(x_2)} + \|P_F(x_1) - P_F(x_2)\|^2 \\
& & \le (2\delta - 1)\|P_F(x_1) - P_F(x_2)\|^2 \\
& & \le 0,
\end{eqnarray*}
provided we choose $\delta \le 1/2$.
\finpf

As before, consider sets $F_1,F_2,\ldots,F_m \subset \E$ and a point $\bx \in \cap_i F_i$,
but now let us suppose moreover that each set $F_i$ is prox-regular at $\bx$.  Define a function
$f \colon \E \rt \R$ by
\bmye \label{mean-squared}
f = \frac{1}{2m} \sum_{i=1}^m d_{F_i}^2.
\emye
This function is half the {\em mean-squared-distance} from the point $x$ to the set system $\{F_i\}$.
According to the preceding result, $f$ is continuously differentiable around $\bx$, and its gradient
\bmye \label{gradient}
\nabla f = \frac{1}{m} \sum_{i=1}^m (I - P_{F_i}) = I - \frac{1}{m} \sum_{i=1}^m P_{F_i}
\emye
is Lipschitz continuous with constant $1$ on a neighborhood of $\bx$.  The method of averaged projections constructs the new iterate $x_+ \in \E$ from the old iterate $x \in \E$ via the update
\bmye \label{e:J4}
  x_+ = \frac{1}{m} \sum_{i=1}^m P_{F_i}(x) = x - \nabla f(x),
\emye
so we can interpret it as the method of steepest descent with a step size of one when the sets $F_i$ are all prox-regular.
To understand its convergence, we return to our strong regularity assumption.

The condition modulus controls the behavior of normal vectors not just at the point $\bx$ but also at nearby points.

\begin{prop}[local effect of condition modulus] \label{nearby}
\hfill
\mbox{Consider closed} sets $F_1,F_2,\ldots,F_m \subset \E$ having strongly regular intersection at a point
$\bx \in \cap F_i$, and any constant
\[
k > \mbox{\rm cond}(F_1,F_2,\ldots,F_m | \bx).
\]
Then for any points $x_i \in F_i$ near $\bx$, any vectors $y_i \in N_{F_i}(x_i)$ (for $i=1,2,\ldots,m$) satisfy the inequality
\[
\sqrt{\sum_i \|y_i\|^2} \le k \Big\| \sum_i y_i \Big\|.
\]
\end{prop}

\pf
If the result fails, then we can find sequences of points $x_i^r \rt \bx$ in $F_i$ and sequences of vectors $y_i^r \in N_{F_i}(x_i)$ (for $i=1,2,\ldots,m$) satisfying the inequality
\[
\sqrt{\sum_i \|y_i^r\|^2} > k \Big\| \sum_i y_i^r \Big\|
\]
for all $r=1,2,\ldots$.  Define new vectors
\[
u_j^r = \frac{1}{\sqrt{\sum_i \|y_i^r\|^2}} \: y_j^r \in N_{F_i}(x_i)
\]
for each index $j=1,2,\ldots,m$ and $r$.  Notice
\[
\sum_i \|u_i^r\|^2 = 1 ~~\mbox{and}~~ \Big\| \sum_i u_i^r \Big\| < \frac{1}{k}.
\]
For each $i=1,2,\ldots,$ the sequence $u_i^1,u_i^2,\ldots$ is bounded, so after taking subsequences we can suppose it converges to some vector $u_i \in \E$, and since the normal cone $N_{F_i}$ is closed as a set-valued mapping from $F_i$ to $\E$, we deduce
$u_i \in N_{F_i}(\bx)$.  But then we have
\[
\sum_i \|u_i\|^2 = 1 ~~\mbox{and}~~ \Big\| \sum_i u_i \Big\| \le \frac{1}{k},
\]
contradicting the definition of the condition modulus $\mbox{cond}(F_1,F_2,\ldots,F_m | \bx)$.  The result follows.
\finpf

The size of the gradient of the mean-squared-distance function $f$, defined by equation (\ref{mean-squared}), is closely related to the value of the function near a point of
strongly regular intersection.  To be precise, we have the following result.

\begin{prop}[gradient of mean-squared-distance] \label{value-gradient}
Consider prox-regular sets $F_1,F_2,\ldots,F_m \subset \E$ having strongly regular intersection at a point $\bx \in \cap F_i$, and any constant
\[
k > \mbox{\rm cond}(F_1,F_2,\ldots,F_m | \bx).
\]Then on a neighborhood of $\bx$, the mean-squared-distance function
\[
f = \frac{1}{2m} \sum_{i=1}^m d_{F_i}^2
\]
satisfies the inequalities
\bmye \label{e:J5}
  \frac{1}{2} \|\nabla f\|^2 \le f \le \frac{k^2m}{2} \|\nabla f\|^2.
\emye
\end{prop}

\pf
Consider any point $x \in \E$ near $\bx$.  By equation (\ref{gradient}), we know
\[
\nabla f(x) = \frac{1}{m} \sum_{i=1}^m y_i,
\]
where
\[
y_i = x_i - P_{F_i}(x_i) \in N_{F_i}(P_{F_i}(x_i))
\]
for each $i=1,2,\ldots,m$.  By definition, we have
\[
f(x) = \frac{1}{2m} \sum_{i=1}^m \|y_i\|^2.
\]
Using inequality (\ref{cauchy}), we obtain
\[
m^2 \|\nabla f(x)\|^2 = \Big\| \sum_{i=1}^m y_i \Big\|^2 \le m \sum_{i=1}^m \| y_i \|^2 = 2m^2 f(x)
\]
On the other hand, since $x$ is near $\bx$, so are the projections $P_{F_i}(x)$, so
\[
2mf(x) = \sum_i \|y_i\|^2 \le k^2 \Big\| \sum_i y_i \Big\|^2 = k^2 m^2 \| \nabla f(x) \|^2.
\]
by Proposition \ref{nearby}.  The result now follows.
\finpf

A standard argument now gives the main result of this section.

\begin{thm}[Q-linear convergence for averaged projections] \label{t:qlin}\hfill \mbox{}
Consider prox-regular sets $F_1,F_2,\ldots,F_m \subset \E$ having strongly regular intersection at a point $\bx \in \cap F_i$, and any constant $k > \mbox{\rm cond}(F_1,F_2,\ldots,F_m | \bx)$.  Then, starting from any point near $\bx$, one iteration of the method of averaged projections reduces the mean-squared-distance
\[
f = \frac{1}{2m} \sum_{i=1}^m d_{F_i}^2
\]
by a factor of at least $1-\frac{1}{k^2 m}$.
\end{thm}

\pf
Consider any point $x \in \E$ near $\bx$.  The function $f$ is continuously differentiable around the minimizer $\bx$, so the gradient $\nabla f(x)$ must be small, and hence the new iterate $x_+ = x - \nabla f(x)$ must also be near $\bx$.  Hence, as we observed after equation (\ref{gradient}), the gradient $\nabla f$ has Lipschitz constant one on a neighborhood of the line segment $[x,x_+]$.  Consequently,
\begin{eqnarray*}
\lefteqn{
f(x_+) - f(x)
} \\
& & = \int_0^1 \frac{d}{dt} f(x-t \nabla f(x))\, dt \\
& & = \int_0^1 \ip{-\nabla f(x)}{\nabla f(x-t \nabla f(x))}\, dt \\
& & = \int_0^1
\Big(
-\|\nabla f(x)\|^2 +
\ip{\nabla f(x)}{\nabla f(x) - \nabla f(x-t \nabla f(x))}
\Big)
\, dt \\
& & \le -\|\nabla f(x)\|^2 + \int_0^1
\|\nabla f(x)\| \cdot \|\nabla f(x) - \nabla f(x-t \nabla f(x))\|
\, dt  \\
& & \le -\|\nabla f(x)\|^2 + \int_0^1 \|\nabla f(x)\|^2 t \, dt \\
& & = -\frac{1}{2}\|\nabla f(x)\|^2 \\
& & \le -\frac{1}{k^2 m} f(x),
\end{eqnarray*}
using Proposition \ref{value-gradient}.
\finpf

A simple induction argument now gives an independent proof in the prox-regular case that the method of averaged projections converges linearly to a point in the intersection of the given sets.  Specifically, the result above shows that mean-squared-distance $f(x_k)$ decreases by at least a constant factor at each iteration, and Proposition \ref{value-gradient} shows that the size of the step $\|\nabla f(x_k)\|$ also decreases by a constant factor.  Hence the sequence $(x_k)$ must converge R-linearly to a point in the intersection.

Comparing this result to Theorem \ref{lin-con} (linear convergence of averaged projections),
we see that the predicted rates of linear convergence are the same.
Theorem \ref{lin-con} guarantees that the squared distance to the intersection
converges to zero with R-linear rate $c^2$ (for any constant $c \in (\bar c,1)$).
The argument gives no guarantee about improvements in a particular iteration:
it only describes the asymptotic behavior of the iterates.
By contrast, the argument of Theorem \ref{t:qlin}, with the added assumption of prox-regularity,
guarantees the same behavior but with the stronger information
that the mean-squared-distance decreases monotonically to zero with Q-linear rate $c^2$.
In particular, each iteration must decrease the mean-squared-distance.

\section{A Numerical Example}\label{examples}

In this final section, we give a numerical illustration showing the linear convergence of
alternating and averaged projections algorithms.
Some major problems in signal or image processing
come down to reconstructing an object
from as few linear measurements as possible.
Several recovery procedures from randomly sampled signals
have been proved to be effective when combined with sparsity constraints
(see for instance the recent developments of compressed sensing
\cite{Cand05},\cite{Don06}).
These optimization problems can be cast as linear programs.  However
for extremely large and/or nonlinear problems, projection methods become attractive
alternatives.   In the spirit of compressive sampling we use projection algorithms to
optimize the compression matrix.  This speculative example is meant simply to illustrate the theory rather than make any claim on real applications.


We consider the decomposition of images $x\in R^n$ as $x=Wz$
where $W \in \R^{n\times m}$ ($n<m$) is a ``dictionary''
(that is, a redundant collection of basis vectors).
Compressed sensing consists in linearly reducing
$x$ to $y=Px=PWz$  with the help of a compression matrix $P\in \R^{d\times n}$
(with $d\ll n$); the inverse operation is to recover $x$ (or $z$) from~$y$.
Compressed sensing theory gives sparsity conditions on $z$
to ensure exact recovery \cite{Cand05},\cite{Don06}. Reference
\cite{Cand05} in fact proposes a recovery algorithm based on alternating projections
(on two convex sets).
In general, we might want to design a specific sensing matrix $P$
adapted to $W$, to ease this recovery process.
An initial investigation on this question is \cite{Ela06};
we suggest here another direction, inspired by \cite{CandesRomberg07}
and \cite{Gab07},
where averaged projections naturally appear.

Candes and Romberg \cite{CandesRomberg07} showed that, under
orthogonality conditions, sparse recovery is more efficient when the entries $|(PW)_{ij}|$ are small.
One could thus use the componentwise $\ell_{\infty}$ norm of $PW$
as a measure of quality of $P$.
This leads to the following feasibility problem:
to find $U=PW$ such that $U\trans{U}=I$ and with the
infinity norm constraint $\norminf{U} \leq \alpha$
(for a fixed tolerance $\alpha$).
The sets corresponding to these constraints are given by
\begin{eqnarray*}
  L&=&\{U\in \R^{d\times m} : U=PW\},\\
  M&=& \{U\in \R^{d\times m} : U\trans{U}=I\},\\
  C&=&\{U\in \R^{d\times m} : \norminf{U}\leq \alpha\}.
\end{eqnarray*}
The first set $L$ is a subspace,
the second set $M$ is a smooth manifold while the third $C$ is convex;
hence the three are prox-regular.
Moreover we can easily compute the projections.
The projection onto the linear subspace $L$ can be computed with a pseudo-inverse.
The manifold $M$ corresponds to the set of matrices $U$ whose singular values are all ones;
it turns out that naturally
the projection onto $M$ is obtained by computing the singular value decomposition of $U$,
and setting singular values to $1$ (apply for example Theorem 27 of \cite{Lew06Alternating}).  Finally the projection onto $C$ comes by shrinking entries of $U$
(specifically, we operate $\min\{\max\{u_{ij},-\alpha\},\alpha\}$ for each entry $u_{ij}$).
This feasibility problem can thus be treated by projection algorithms, and hopefully
a matrix $U\in L\cap M\cap C$ will correspond to a good compression matrix $P$.

\begin{figure}\label{f:OCS}
  \begin{center}
     \includegraphics[width=12cm, height=7cm]{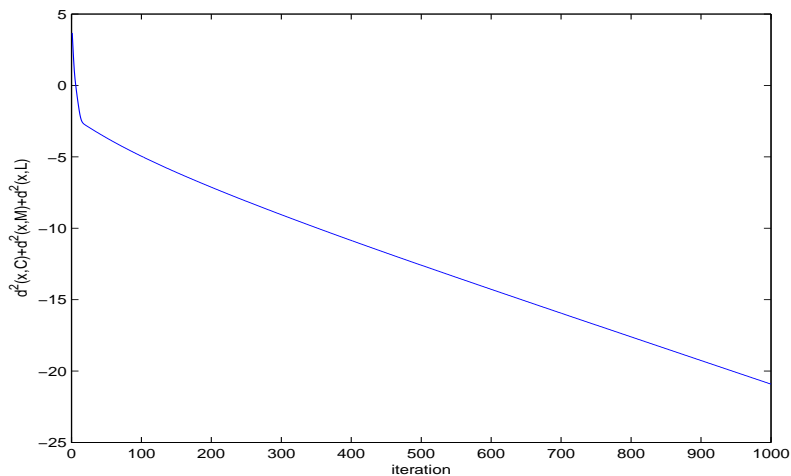}
  \end{center}
  \caption{Convergence of averaged projection algorithm
    for designing compression matrix in compressed sensing.}
\end{figure}

To illustrate this problem, we generate random entries (normally distributed)
of the dictionary $W$ (size $128\times 512$, redundancy factor $4$)
and of an initial iterate $U_0\in L$. We fix $\alpha=0.1$,
and we run the averaged projection algorithm
which computes a sequence of $U_k$.
Figure \ref{f:OCS} shows
\[
\log_{10} \ f(U_k) \qquad \text{with }f(U)=\frac{1}{6}(d^2_{L}(U)+d^2_{M}(U)+d^2_{C}(U))
\]
for each iteration $k$. We also observe that the ratio
\[
f(U_{k+1})/f(U_{k+1}) < 0.9627
\]
for all iterations $k$, showing the expected Q-linear convergence.
We note that working on random test cases is of interest
for our simple testing of averaged projections:
though we cannot guarantee in fact that the intersection of the three sets is strongly regular,
randomness seems to prevent irregular solutions, providing $\alpha$ is not too small.
So in this situation, it is likely that the algorithm will converge linearly;
this is indeed what we observe in Figure~\ref{f:OCS}.
We note furthermore that we tested alternating projections on this problem
(involving three sets, so not explicitly covered by Theorem \ref{key}).  We observed that the method is still converging linearly in practice, and again, the rate is better than for averaged projections.

This example illustrates how the projection algorithm behaves
on random feasibility problems of this type.
However the potential benefits of using optimized compression matrix versus
random compression matrix in practice are still unclear.
Further study and more complete testing have to be done for these questions;
this is beyond the scope of this paper.

\end{document}